\documentstyle[11pt,A4]{article}
\begin{document}
\newtheorem{algorithm}{Algorithm}
\newtheorem{lemma}{Lemma}
\newtheorem{theorem}{Theorem}
\newtheorem{corollary}{Corollary}
\newtheorem{definition}{Definition}
\newcommand{\beq}{\begin{equation}}
\newcommand{\eeq}{\end{equation}}
\newcommand{\QED}{{\hfill$\Box$}}
\newcommand{\E}{{$'$}}
\newcommand{\veps}{\varepsilon}
\newcommand{\ou}{{\overline{u}}}
\newcommand{\oy}{{\overline{y}}}
\newcommand{\ozz}{{\overline{z}}}
\newcommand{\cd}{{\tilde d}}
\newcommand{\cuu}{{\tilde u}}
\newcommand{\cx}{{\tilde x}}
\newcommand{\cy}{{\tilde y}}
\newcommand{\coz}{{\tilde{\overline{z}}}}
\newcommand{\cQ}{{\tilde Q}}
\newcommand{\cR}{{\tilde R}}
\newcommand{\ccR}{{\tilde R'}}
\newcommand{\cRb}{{\tilde R_b}}
\newcommand{\cRt}{{\tilde R_t}}
\newcommand{\cU}{{\tilde U}}
\newcommand{\hA}{{\widehat A}}
\newcommand{\hQ}{{\widehat Q}}
\newcommand{\hb}{{\widehat b}}
\newcommand{\Al}{{A_{-1}}}
\newcommand{\AlT}{{A_{-1}^T}}
\newcommand{\Rbb}{{R_b}}
\newcommand{\Rtt}{{R_t}}
\newcommand{\rll}{{r_{1,1}}}
\newcommand{\rls}{{r_{1,1}^2}}
\newcommand{\calA}{{\cal A}}
\newcommand{\calB}{{\cal B}}
\newcommand{\calF}{{\cal F}}
\newcommand{\disp}{{\cal D}}
\newcommand{\myhline}{\hline&\\[-12pt]}	%
\bibliographystyle{plain}

\title{A Weakly Stable Algorithm for General Toeplitz Systems%
\thanks{Copyright \copyright\ 1993--2010 the authors.
The results were announced in
{\em Parallel algorithms and numerical stability for Toeplitz
systems}, invited paper presented by Brent at the
SIAM Conference on Linear Algebra in Signals,
Systems and Control, Seattle, August 16-19, 1993.
		\hspace*{\fill} \hbox{rpb143tr typeset using \LaTeX}}}

\author{Adam W.\ Bojanczyk\\
	School of Electrical Engineering\\
	Cornell University\\
	Ithaca, NY 14853-5401
\and
	Richard P.\ Brent\\
	Computer Sciences Laboratory\\
	Australian National University\\
	Canberra, ACT 0200
\and
	Frank R.\ de Hoog\\
	Division of Mathematics and Statistics\\
	CSIRO, GPO Box 1965\\
	Canberra, ACT 2601\\[15pt]}

\date{Report TR-CS-93-15\\
	6 August 1993\\
	(revised 24 June 1994)}

\maketitle
\vspace{-15pt}
\thispagestyle{empty}			%

\begin{abstract}

We show that a fast algorithm for the
$QR$ factorization of a Toeplitz or Hankel matrix $A$ is weakly stable in the
sense that $R^TR$ is close to $A^TA$. Thus, when the algorithm is used
to solve the semi-normal equations $R^TRx = A^Tb$,
we obtain a weakly stable method for the solution
of a nonsingular Toeplitz or Hankel linear system $Ax = b$.
The algorithm also applies to the solution of the full-rank
Toeplitz or Hankel least squares problem $\min \| Ax-b \|_2$.
\smallskip

\smallskip


\end{abstract}

\section{Introduction}
\label{sec:Intro}

Toeplitz linear systems arise in many applications, and there
are many algorithms which solve nonsingular $n \times n$ Toeplitz systems
\[	Ax = b	\]
in $O(n^2)$ arithmetic operations~\cite{Bareiss69,rpb111,%
Heinig-JR,Heinig-Rost,Kailath86,Kolmogorov41,Levinson47,Rissanen74,%
Schur17,Szego39,Trench64,Watson73,Wiener49,Zohar69,Zohar74}.
Some algorithms are restricted to
symmetric systems ($A = A^T$) and others apply to general Toeplitz
systems.  Because of their recursive nature, most $O(n^2)$ algorithms
assume that all leading principal submatrices of $A$ are nonsingular,
and break down if this is not the case.
These algorithms are generally unstable,
because a leading principal submatrix
may be poorly conditioned even if $A$ is well conditioned.
Thus, stability results often depend on the assumption that
$A$ is symmetric positive definite, in which case the leading
principal submatrices are at least as well conditioned as $A$.

Asymptotically faster algorithms
exist~\cite{Ammar88,Bitmead80,rpb059,Chun91,%
Gutknecht92,Hoog87,Morf80,Musicus84}.
Sometimes these algorithms are called
{\em superfast}~\cite{Ammar88}.
We avoid this terminology because, even though the algorithms require
only $O(n(\log n)^2)$ arithmetic operations,
they may be slower than
$O(n^2)$ algorithms for $n < 256$ (see~\cite{Ammar88,Chun91,Sexton82}).
We prefer the term {\em asymptotically fast}.

The numerical stability properties of asymptotically fast algorithms
are generally either bad~\cite{Bunch85} or unknown,
although some positive partial results have been obtained
recently~\cite{Gutknecht93b}.
Attempts to stabilise asymptotically fast algorithms by look-ahead techniques
have been made~\cite{Gutknecht92}, but the look-ahead algorithms are
complicated and their worst-case behaviour is unclear.
Thus, we do not consider asymptotically fast algorithms further,
but restrict our attention to $O(n^2)$ algorithms.

We are concerned with direct methods for the general case,
where $A$ is any nonsingular Toeplitz matrix.
In this case no $O(n^2)$ algorithm has been proved to be stable.
For example, the algorithm of Bareiss~\cite{Bareiss69} has stability
properties similar to those of Gaussian elimination
without pivoting~\cite{rpb126,rpb144,Sweet82,Sweet93},
so is unstable and breaks down if a leading principal minor vanishes.
Several authors have suggested the introduction of
pivoting or look-ahead (with block steps) in the Bareiss and Levinson
algorithms~\cite{Chan92a,Chan92b,Freund93c,Freund93a,Freund93b,Sweet90,Sweet93},
and this is often successful in practice,
but in the worst case the overhead is $O(n^3)$ operations.
The recent algorithm GKO of Gohberg, Kailath and Olshevsky~\cite{Gohberg94}
may be as stable as Gaussian elimination with partial pivoting,
but an error analysis has
not been published. %

In an attempt to achieve stability without pivoting or look-ahead,
it is natural to consider algorithms for computing an
orthogonal factorization
\beq
	A = QR						\label{eq:AQR}
\eeq
of $A$. The first such $O(n^2)$ algorithm was introduced by
Sweet~\cite{Sweet82,Sweet84}. Unfortunately, Sweet's algorithm
depends on the condition of a certain submatrices of $A$,
so is unstable~\cite{rpb092,Luk87}.
Other $O(n^2)$ algorithms for computing the matrices $Q$ and $R$ or $R^{-1}$
in~(\ref{eq:AQR}) were given by
Bojanczyk, Brent and de~Hoog~\cite{rpb092},
Chun, Kailath and Lev-Ari~\cite{Chun87},
Cybenko~\cite{Cybenko87},
Luk and Qiao~\cite{Luk87,Qiao88},
and Nagy~\cite{Nagy93}.
To our knowledge none of them has
been shown to be stable.
In several cases examples show that they are not stable.
Unlike the classical $O(n^3)$ Givens or Householder algorithms,
the $O(n^2)$ algorithms do not form $Q$
in a numerically stable manner
as a product of matrices which are (close to) orthogonal.

Numerical experiments with the algorithm of Bojanczyk, Brent and de~Hoog
(BBH for short) suggest that the
cause of instability is the method for computing
the orthogonal matrix $Q$;
the computed upper triangular matrix $\cR$
is about as good as can be obtained by performing a Cholesky
factorization of $A^TA$, provided the downdates involved in the algorithm are
implemented in a certain way (see~\S\ref{sec:Downdating}).
This result is proved in \S\ref{sec:Main}.
As a consequence, in \S\ref{sec:sneq} we show how the method of
semi-normal equations (i.e.~the solution of $R^TRx = A^Tb$) can be used
to give a weakly stable algorithm for the solution of general Toeplitz
or Hankel systems.
The result also applies to the solution of full-rank Toeplitz or Hankel least
squares problems. For a discussion of the rank-deficient case and a
``look-ahead'' modification of Cybenko's algorithm, see~\cite{Hansen93}.

In \S\ref{sec:Notation} we introduce some notation and conventions.
The concepts of stability and weak stability are defined in
\S\ref{sec:Stability}.
The Cholesky downdating problem, which arises in the BBH algorithm,
is discussed in \S\ref{sec:Downdating}.
Numerical results are discussed in \S\ref{sec:Numer},
and some conclusions are given in \S\ref{sec:Conc}.

If $H$ is a Hankel matrix, and $J$ is the permutation matrix which
on premultiplication reverses the order of the rows of a matrix,
then $JH$ is Toeplitz. Also, $(JH)^T(JH) = H^TH$.
Thus, our results apply equally to Hankel and Toeplitz matrices.
Our results might also be extended to more general classes of matrices
with a displacement
structure~\cite{Kailath94,Kailath79,Kailath-Sayed,Kailath78}.
For simplicity we restrict our attention to the Toeplitz case.

\pagebreak[4]
\section{Notation and Conventions}
\label{sec:Notation}

Let

\[
A = \left(\begin{array}{ccc}
a_0	& \cdots	& a_{n-1}	\\
\vdots	& \ddots	& \vdots	\\
a_{1-m}	& \cdots	& a_{n-m}
\end{array}\right)
\]

\noindent be a real $m \times n$ Toeplitz matrix, so
$a_{i,j} = a_{j-i}$ for $1 \le i \le m$, $1 \le j \le n$.
We assume that $m \ge n$ and that $A$ has rank $n$.
Thus $A^TA$ has a Cholesky factorization
$A^TA = R^TR$,
where $R$ is an upper triangular $n \times n$ matrix.
We assume that the diagonal elements of $R$ are positive, so
$R$ is unique. Also,
$A = QR$,
where $Q$ is an $m \times n$ matrix with orthonormal columns.
\vfil

If the singular values of $A$ are $\sigma_1, \ldots, \sigma_n$, where
$\sigma_1 \ge \ldots \ge \sigma_n > 0$,
then the spectral condition number of $A$ is
\[	\kappa = \kappa_2(A) = \sigma_1/\sigma_n .\]
For convenience in stating the error bounds, we often assume that $\sigma_1$
is of order unity, which can be achieved by scaling.
\vfil

A {\em displacement operator}
$\disp : {\Re}^{m\times n} \to {\Re}^{(m-1)\times(n-1)}$
is defined as follows:
for any $m \times n$ matrix $B$, $\disp(B) = C$, where $C$ is
the $(m-1) \times (n-1)$ matrix with entries
$c_{i,j} = b_{i+1,j+1} - b_{i,j}$,
$1 \le i < m, 1 \le j < n$.
Note that $\disp B = 0$ iff $B$ is Toeplitz.
\vfil

Let $\veps$ be the machine precision. In our analysis we neglect
terms of order $O(\veps^2)$. This can be justified by considering
our bounds as asymptotic expansions in the (sufficiently small) parameter
$\veps$.
\vfil

It is often convenient to
subsume a polynomial in $m$ and/or $n$ into the ``$O$'' notation,
and indicate this by a subscript.
Thus, an error bound of the form
\[	\|E\| = O_m(\veps)	\]
 means that
\[	\| E \| \le P(m)\veps	\]
for some polynomial $P$ and all sufficiently small $\veps$.
This notation is useful because minor changes in the algorithm
or changes in the choice of norm will be absorbed by a change in
the polynomial $P(m)$.
It is often stated (e.g.~by Bj\"orck~\cite{Bjorck91})
that {\em the primary purpose of rounding error analysis is insight},
and insight can be aided by the suppression of superfluous details.%
\footnote{Wilkinson~\cite{Wilkinson61} states	%
``$\ldots$ there is a danger that the essential simplicity of the
error analysis may be obscured by an excess of detail.''}
\vfil

If the error bound depends on $\kappa$ then this will
be mentioned explicitly (e.g.~$\| E \| = O_m(\kappa\veps)$).
The meaning of ``sufficiently small'' may depend on $\kappa$
For example, we may need $\veps < 1/\kappa^2$.
\smallskip
\vfil

\noindent We distinguish several classes of numerical quantities~--
\begin{enumerate}
\item Exact values, e.g.~input data such as~$a_i$.
\item Computed values, usually indicated by a tilde, e.g.~$\tilde u_i$.
\item Perturbed values given by error analysis, usually indicated by
a hat, e.g.~${\widehat a}_{i,j}$, or by
the argument $\veps$, e.g.~$a_{i,j}(\veps)$.
These are not computed, but the error analysis shows that
they exist and gives bounds on their difference from the corresponding
exact values.
Sometimes the perturbations are of computed quantities,
e.g.~$\tilde u_i(\veps)$.
\end{enumerate}
\vfil

\pagebreak[4]
\section{Stability and Weak Stability}
\label{sec:Stability}

In this section we give definitions of stability and weak stability
of algorithms for solving linear systems.

Consider algorithms for solving a nonsingular, $n \times n$ linear
system $Ax = b$, so $m = n$.
There are many definitions of %
numerical stability in the literature, for example
\cite{Bjorck87,Bjorck91,rpb126,Bojanczyk91,Bunch85,Cybenko80,%
Golub89,Jankowski77,Miller80,Paige73,Stewart73}.
Definitions~\ref{def:stable} and~\ref{def:weak} below are taken from
Bunch~\cite{Bunch87}. %

\begin{definition}
\label{def:stable}
An algorithm for solving linear equations
is {\em stable} for a class of matrices $\calA$ if for each $A$ in
$\calA$ and for each $b$ the computed solution $\cx$ to $Ax = b$
satisfies $\hA\cx = \hb$, where $\hA$ is close to $A$ and $\hb$ is
close to $b$.
\end{definition}

Definition~\ref{def:stable} says that, for stability,
the {\em computed} solution
has to be the {\em exact} solution of a problem
which is close to the original problem.
This is the classical {\em backward stability} of
Wilkinson~\cite{Wilkinson61,Wilkinson63,Wilkinson65}.
We interpret ``close'' to mean close
in the relative sense in some norm, i.e.
\beq
	\|\hA - A\|/\|A\| = O_n(\veps),\;
	\|\hb - b\|/\|b\| = O_n(\veps).			\label{eq:close}
\eeq

Note that the matrix $\hA$ is not required to be in the class $\calA$.
For example, $\calA$ might be the class of nonsingular Toeplitz matrices,
but $\hA$ is not required to be a Toeplitz matrix. If we require
$\hA \in \calA$ we get what Bunch~\cite{Bunch87} calls {\em strong stability}.
For a discussion of the difference between stability and strong stability
for Toeplitz algorithms, see~\cite{Higham92,Varah92}.

Stability does not imply that the computed
solution $\cx$ is close to the exact solution $x$, unless the problem
is well-conditioned.
Provided $\kappa\veps$ is sufficiently small, stability implies that
\beq
\|\cx - x\|/\|x\| = O_n(\kappa\veps).			\label{eq:relerr}
\eeq
For more precise results,
see~Bunch~\cite{Bunch87} and~Wilkinson~\cite{Wilkinson61}.

As an example, consider the method of Gaussian elimination.
Wilkinson~\cite{Wilkinson61} shows that
\beq
\|\hA - A\|/\|A\| = O_n(g\veps),			\label{eq:GEerr}
\eeq
where $g = g(n)$ is the ``growth factor''.
$g$ depends on whether partial or complete pivoting is used.
In practice $g$ is usually moderate, even for partial pivoting.
However, a well-known example shows that $g(n) = 2^{n-1}$ is possible
for partial pivoting, and recently
it has been shown that examples where $g(n)$ grows exponentially
with $n$ may arise in applications, e.g.~for linear systems arising from
boundary value problems.	%
Even for complete pivoting,
it has not been {\em proved} that $g(n)$ is bounded
by a polynomial in $n$.
Wilkinson~\cite{Wilkinson61}	%
showed that $g(n) \le n^{(\log n)/4 + O(1)}$,
and Gould~\cite{Gould91} recently showed that $g(n) > n$ is possible
for $n > 12$; there is still a large gap between these results.
Thus, to be sure that Gaussian elimination satisfies
Definition~\ref{def:stable}, we must restrict $\calA$
to the class of matrices for which $g$ is $O_n(1)$.
In practice this is not a problem, because $g$ can easily be checked
{\em a~posteriori}.

Although stability is desirable, it is more than we can prove for many
useful algorithms. Thus, following Bunch~\cite{Bunch87}, we define the
(weaker, but still useful)
property of {\em weak stability}.

\begin{definition}
\label{def:weak}
An algorithm for solving linear equations
is {\em weakly stable} for a class of matrices $\calA$ if for each
well-conditioned $A$ in
$\calA$ and for each $b$ the computed solution $\cx$ to $Ax = b$
is such that $\|\cx-x\|/\|x\|$ is small.
\end{definition}

In Definition~\ref{def:weak}, we take ``small'' to mean $O_n(\veps)$,
and ``well-conditioned'' to mean that $\kappa(A)$ is $O_n(1)$,
i.e.~is bounded by a polynomial in $n$.
From~(\ref{eq:relerr}), stability implies weak stability.

Define the {\em residual} $r = A\cx - b$.
It is well-known~\cite{Wilkinson63} that
\beq
{1 \over \kappa}{\|r\|\over\|b\|} \le
{\|\cx-x\| \over \|x\|} \le \kappa{\|r\|\over\|b\|} .
\eeq
Thus, for well-conditioned $A$, $\|\cx-x\|/\|x\|$ is small if and only if
$\|r\|/\|b\|$ is small. This observation clearly
leads to an alternative definition of weak stability:

\begin{definition}
\label{def:weak2}
An algorithm for solving linear equations
is {\em weakly stable} for a class of matrices $\calA$ if for each
well-conditioned $A$ in
$\calA$ and for each $b$ the computed solution $\cx$ to $Ax = b$
is such that $\|A\cx - b\|/\|b\|$ is small.
\end{definition}

Now consider computation of the Cholesky factor $R$ of $A^TA$, where
$A$ is an $m \times n$ matrix of full rank $n$.
A good $O(mn^2)$ algorithm is
to compute the $QR$ factorization
\[	A = QR	\]
of $A$ using Householder or Givens transformations~\cite{Golub89}.
It can be shown~\cite{Wilkinson63}
that the computed matrices $\cQ$, $\cR$ satisfy
\beq
	\hA = \hQ\cR				\label{eq:QRerror}
\eeq
where $\hQ^T\hQ = I$,
$\cQ$ is close to $\hQ$, and $\hA$ is close to $A$.
Thus, the algorithm is stable in the
sense of backward error analysis.
Note that $\|A^TA - \cR^T\cR\|/\|A^TA\|$ is small,
but $\|\cQ - Q\|$ and $\|\cR - R\|/\|R\|$ are not necessarily small.
Bounds on $\|\cQ - Q\|$ and $\|\cR - R\|/\|R\|$
depend on $\kappa$, and
are discussed in~\cite{Golub65,Stewart77,Wilkinson65}.

A different algorithm is to compute (the upper triangular part of) $A^TA$,
and then compute the Cholesky factorization of $A^TA$ by the usual
(stable) algorithm. The computed result $\cR$ is such that
$\cR^T\cR$ is close to $A^TA$. However, this does not
imply the existence of $\hA$ and $\hQ$ such that~(\ref{eq:QRerror}) holds
(with $\hA$ close to $A$ and some $\hQ$ with $\hQ^T\hQ = I$)
unless $A$ is well-conditioned~\cite{Stewart79}. By analogy with
Definition~\ref{def:weak2} above,
we may say that Cholesky factorization of $A^TA$ gives a
{\em weakly stable} algorithm for computing $R$,
because the ``residual'' $A^TA - \cR^T\cR$ is small.

\section{Cholesky Updating and Downdating}
\label{sec:Downdating}

\subsection{Updating}
\label{subsec:Updating}

The Cholesky updating problem is: given an upper triangular matrix
$R \in \Re^{n \times n}$ and a vector $x \in \Re^n$,
find an upper triangular matrix $U$ such that
\beq
	U^TU = R^TR + xx^T.                  		\label{eq:update}
\eeq
The updating problem can be solved in a numerically stable manner
by transforming the matrix
$\left(\begin{array}{c}
x^T \\
R
\end{array}\right)$
to upper triangular form
$\left(\begin{array}{c}
U \\
0^T
\end{array}\right)$
by applying a sequence of plane rotations on the left.
For details, see~\cite{Golub89}.

\subsection{Downdating}
\label{subsec:Downdating}

The Cholesky downdating problem is: given a upper triangular matrix
$R \in \Re^{n \times n}$ and a vector $x \in \Re^n$
such that $R^TR - xx^T$ is positive definite,
find an upper triangular matrix $U$ such that
\beq
	U^TU = R^TR - xx^T.				\label{eq:downdate}
\eeq
Proceeding formally, we can obtain~(\ref{eq:downdate}) from~(\ref{eq:update})
by replacing $x$ by $ix$. However, the numerical properties of the
updating and downdating problems are very different.
The condition that $R^TR - xx^T$
be positive semi-definite is necessary for the existence of
a real $U$ satisfying~(\ref{eq:downdate}).
Thus, we would expect the downdating problem to be illconditioned
if $R^TR - xx^T$ has small singular values,
and Stewart~\cite{Stewart79} shows that this is true.

There are several algorithms for the Cholesky downdating
problem~\cite{Bischof93,rpb095,Bojanczyk88,Bojanczyk91,Dongarra78,%
Fletcher74,Gill74,Golub89,Pan90,Saunders72} and we shall not discuss
them in detail here. What is relevant to us is the error analysis.
To simplify the statement of the error bounds, suppose that $\|R\| = O_n(1)$,
which implies that $\|x\| = O_n(1)$.
Observe that, in exact arithmetic,
there is an orthogonal matrix $Q$ such that
\beq
\left(\begin{array}{c}
x^T \\
U
\end{array}\right) =
Q \left(\begin{array}{c}
R \\
0^T
\end{array}\right).
\eeq
Suppose the computed upper triangular matrix is $\cU$.
Stewart~\cite{Stewart79} shows that,
for the ``Linpack'' algorithm~\cite{Dongarra78},

\beq
\left(\begin{array}{c}
x^T(\veps) \\
\cU(\veps)
\end{array}\right) =
Q(\veps) \left(\begin{array}{c}
R \\
0^T
\end{array}\right),					\label{eq:linperr}
\eeq
where $Q(\veps)$ is an exactly orthogonal matrix,
\beq
\|x(\veps) - x\| = O_n(\veps),			\label{eq:linperr1}
\eeq
and
\beq
\|\cU(\veps) - \cU\| = O_n(\veps).			\label{eq:linperr2}
\eeq

We can regard $x(\veps)$ as a (backward) perturbation of the input
data $x$, and $\cU(\veps)$ as a (forward) perturbation of the
computed result $\cU$. Because of this mixture of forward and backward
perturbations, a result of this form is sometimes called a
``mixed'' stability result.
If the problem is illconditioned, the backward perturbation is
more significant than the forward perturbation.

It is important to note that the error analysis does {\em not} show that
the computed matrix $\cU$ is close to the exact result $U$,
or that $Q(\veps)$ is close to $Q$,
unless $U$ is well-conditioned.

A stability result of the same form as~(\ref{eq:linperr}--\ref{eq:linperr2})
has been established by
Bojanczyk, Brent, Van~Dooren and de~Hoog~\cite{rpb095}
for their ``Algorithm~C''.
Because Algorithm~C computes $\cU$ one row at
a time, several (updates and/or) downdates can be pipelined,
which is not the case for the Linpack algorithm%
\footnote{The recent algorithm of Bischof et al~\cite{Bischof93} can
also be pipelined, but we do not know if it gives results which satisfy
equations~(\ref{eq:linperr}--\ref{eq:linperr2})
or~(\ref{eq:downerr}--\ref{eq:Gerr1}).}.

The result (\ref{eq:linperr}--\ref{eq:linperr2}) implies that
\beq
	\cU^T\cU = R^TR - xx^T + G(\veps),	\label{eq:downerr}
\eeq
where
\beq
\|G(\veps)\| = O_n(\veps).				\label{eq:Gerr1}
\eeq
Clearly a similar result holds if a sequence of (updates and) downdates
is performed, provided that the final and intermediate results
are positive definite.

\section{A Weak Stability Result for the BBH Algorithm}
\label{sec:Main}

Our main result (Theorem~\ref{thm:main} below)
is that the BBH algorithm computes $R$ about as well
as would be expected if the Cholesky factorization of $A^TA$ were
computed by the usual (stable) algorithm.
More precisely, the computed $\cR$ satisfies
\beq
A^TA - \cR^T\cR = O_m(\veps\|A^TA\|).			\label{eq:MR}
\eeq
To avoid having to include $\|A^TA\|$ in the error bounds,
we assume for the time being that $\sigma_1(A) = O(1)$.

First consider exact arithmetic.
We partition $A$ in two ways:
\beq
A = \left(\begin{array}{c|c}
a_0	& y^T	\\ \hline
z	& \Al \end{array}\right) =
\left(\begin{array}{c|c}
\Al	& \oy	\\ \myhline
\ozz^T	& a_{n-m} \end{array}\right),			\label{eq:Apart}
\eeq
where $\Al$ is an $(m-1) \times (n-1)$ Toeplitz matrix,
\begin{eqnarray*}
y 	&=& (a_1, \ldots, a_{n-1})^T,\\
z       &=& (a_{-1}, \ldots, a_{1-m})^T,\\
\oy	&=& (a_{n-1}, \ldots, a_{n-m+1})^T,\\
\ozz	&=& (a_{1-m}, \ldots, a_{n-m-1})^T.
\end{eqnarray*}

Similarly, we partition $R$ in two ways:
\beq
R = \left(\begin{array}{c|c}
\rll	& u^T	\\ \hline
0	& \Rbb \end{array}\right) =
\left(\begin{array}{c|c}
\Rtt	& \ou	\\ \myhline
0^T	& r_{n,n} \end{array}\right),			\label{eq:Rpart}
\eeq
where $\Rbb$ and $\Rtt$ are $(n-1) \times (n-1)$ upper triangular matrices,
\begin{eqnarray*}
u	&=& (r_{1,2}, \ldots, r_{1,n})^T,\\
\ou	&=& (r_{1,n}, \ldots, r_{n-1,n})^T.
\end{eqnarray*}

From~(\ref{eq:Apart}),	%
\beq
A^TA = \left(\begin{array}{c|c}
a_0^2 + z^Tz	& a_0y^T + z^T\Al \\ \myhline
a_0y + \AlT z	& \AlT\Al + yy^T
\end{array}\right) =
\left(\begin{array}{c|c}
\AlT\Al + \ozz\,\ozz^T	& \vdots \\ \hline
\cdots			& \cdot
\end{array}\right),					\label{eq:ATA}
\eeq
where the dots indicate entries which are irrelevant for our purposes.
Similarly, from~(\ref{eq:Rpart}),
\beq
R^TR =  \left(\begin{array}{c|c}
\rls	& \rll u^T \\ \myhline
\rll u	& \Rbb^T\Rbb + uu^T
\end{array}\right) =
\left(\begin{array}{c|c}
\Rtt^T\Rtt	& \vdots \\ \hline
\cdots		& \cdot
\end{array}\right).					\label{eq:RTR}
\eeq
Equating $A^TA$ and $R^TR$,
we obtain
\beq
\rls = a_0^2 + z^Tz,					\label{eq:r11}
\eeq
\beq
\rll u = a_0y + \AlT z,					\label{eq:row1}
\eeq
\beq
\AlT\Al + yy^T = \Rbb^T\Rbb + uu^T, 			\label{eq:Rb}
\eeq
and
\beq
\AlT\Al + \ozz\,\ozz^T = \Rtt^T\Rtt.			\label{eq:Rt}
\eeq
Eliminating $\AlT\Al$ from~(\ref{eq:Rb}--\ref{eq:Rt}) gives the relation
\beq
	\Rbb^T\Rbb = \Rtt^T\Rtt + yy^T - uu^T - \ozz\,\ozz^T	\label{eq:r3}
\eeq
which is the basis for the BBH algorithm.
If $\Rtt$ were known, then $\Rbb$ could be computed from~(\ref{eq:r3})
using one Cholesky updating step and two Cholesky downdating steps,
as discussed in \S\ref{sec:Downdating}.
Also, since updating and downdating algorithms can proceed by rows,
knowledge of the first $k$ rows of $\Rtt$ is sufficient to allow
the computation of the first $k$ rows of $\Rbb$.
It is easy to compute the first row of $R$ from
relations~(\ref{eq:r11}--\ref{eq:row1}).
(For future reference, suppose that the computed first row of $R$
is $(\tilde{r}_{1,1}, \cuu^T)$.)
It is clear from~(\ref{eq:Rpart}) that the $k$-th row of $\Rbb$
defines the $(k+1)$-th row of $\Rtt$.
Thus, we can compute $\Rtt$ and $\Rbb$ row by row.
For details see~\cite{rpb092}.

A straightforward extension of the
result~(\ref{eq:downerr}--\ref{eq:Gerr1}) applies
to our problem of computing $\Rtt$ and $\Rbb$.
Provided the ``Algorithm~C'' variant of downdating~\cite{rpb095} is used,
the computed results $\cRt$ and $\cRb$ satisfy
\beq
	\cRb^T\cRb = \cRt^T\cRt + yy^T - \cuu\cuu^T - \ozz\,\ozz^T
			+ G(\veps)
							\label{eq:cr3}
\eeq
where
\beq
\|G(\veps)\| = O_m(\veps).				\label{eq:Gerr}
\eeq
Here $y$, $\ozz$ and $\cuu$ are inputs to the
up/downdating procedures
($\cuu$ may differ slightly from the exact $u$
because $\cuu$ has been computed from~(\ref{eq:r11}--\ref{eq:row1})).
At this point we make no claims about the size of
$\|\cRb - \Rbb\|$ and $\|\cRt - \Rtt\|$. All we need
is that $\cRb$ and $\cRt$ exist
and are bounded for sufficiently small $\veps$.

Because of the algorithm for their computation,
the computed matrices $\cRt$ and $\cRb$ are related
so that we can define the ``computed $R$'', say $\cR$,
in a consistent manner by
\beq
\cR = \left(\begin{array}{c|c}
\tilde{r}_{1,1}	& \cuu^T	\\ \myhline
0	& \cRb \end{array}\right) =
\left(\begin{array}{c|c}
\cRt	& \vdots	\\ \myhline
0^T	& \cdot \end{array}\right).			\label{eq:Rpartc}
\eeq
From~(\ref{eq:Rpartc}) we have
\beq
\cR^T\cR = \left(\begin{array}{c|c}
\tilde{r}_{1,1}^2	& \tilde{r}_{1,1}\cuu^T	\\ \hline
\vdots			& \cRb^T\cRb + \cuu\cuu^T \end{array}\right) =
\left(\begin{array}{c|c}
\cRt^T\cRt	& \vdots	\\ \hline
\cdots		& \cdot \end{array}\right).		\label{eq:RTRe}
\eeq
Recall our definition of the operator $\disp$ in \S\ref{sec:Notation}.
From~(\ref{eq:RTRe}) we have
\beq
\disp(\cR^T\cR) = \cRb^T\cRb + \cuu\cuu^T - \cRt^T\cRt .
\eeq
Thus, from~(\ref{eq:cr3}),
\beq
\disp(\cR^T\cR) = yy^T - \ozz\ozz^T + G(\veps).	\label{eq:DR}
\eeq
Also, from~(\ref{eq:ATA}),
\beq
\disp(A^TA) = yy^T - \ozz\,\ozz^T.			\label{eq:DA}
\eeq
If $E = \cR^T\cR - A^TA$ and $F = \disp(E)$ then,
from~(\ref{eq:DR}--\ref{eq:DA}),
\beq
    F = G(\veps).
\eeq
If $1 \le j \le i \le n$ then,	%
by the definition of $\disp(E)$,
\[
e_{i,j} - e_{i-j+1,1}
	= \sum_{k=1}^{j-1} (e_{i-k+1,j-k+1} - e_{i-k,j-k})
	= f_{i-1,j-1} + f_{i-2,j-2} + \cdots + f_{i-j+1,1}.
\]
The first row of $\cR^T\cR$ is $\tilde{r}_{1,1}(\tilde{r}_{1,1}, \cuu^T)$,
which is close to $r_{1,1}(r_{1,1}, u^T)$,
so the first row of $E$ has norm $O_m(\veps)$.
Also, $E$ is symmetric. It follows that
\[ \| E \| \le (n-1)\| F \| + O_m(\veps) = O_m(\veps). \]
Thus, after scaling to remove our assumption that
$\sigma_1 = O(1)$,
we have proved:

\begin{theorem}
\label{thm:main}
If the BBH algorithm is used with the downdating steps performed as in\\
``Algorithm~C'' of~{\rm{\cite{rpb095}}}, then the computed Cholesky
factor $\cR$ of $A^TA$ satisfies
\beq
\| \cR^T\cR - A^TA \| = O_m(\veps\|A^TA\|).		\label{eq:mainbound}
\eeq
\end{theorem}

Since $A^TA = R^TR$, a comparison of the partitioned forms~(\ref{eq:RTR})
and~(\ref{eq:RTRe}) shows that
\beq
\|\cRb^T\cRb - \Rbb^T\Rbb\| = O_m(\veps\|R^TR\|)	\label{eq:RTR2}
\eeq
and
\beq
\|\cRt^T\cRt - \Rtt^T\Rtt\| = O_m(\veps\|R^TR\|).	\label{eq:RTR3}
\eeq
From~(\ref{eq:mainbound}) and Stewart's perturbation analysis~\cite{Stewart79},
it follows that
\beq
\| \cR - R \|/\|R\| = O_m(\kappa \veps).		\label{eq:Rerr}
\eeq
Note that the condition number $\kappa$ appears in~(\ref{eq:Rerr})
but not in~(\ref{eq:mainbound}--\ref{eq:RTR3}).

\section{The Semi-Normal Equations}
\label{sec:sneq}

Suppose that our aim is to solve a nonsingular $n \times n$
Toeplitz linear system
\beq
Ax = b,                                        		\label{eq:Axb}
\eeq
using $O(n^2)$ arithmetic operations.
In exact arithmetic, the {\em normal equations}
\beq
A^TAx = A^Tb                     			\label{eq:normal}
\eeq
and the {\em semi-normal equations}
\beq
R^TRx = A^Tb,						\label{eq:seminormal}
\eeq
where $R$ satisfies~(\ref{eq:AQR}),
are equivalent to~(\ref{eq:Axb}).

In most circumstances the use of the normal or semi-normal equations is
not recommended, because the condition number $\kappa(A^TA)$
may be as large as $\kappa(A)^2$
(see {\S}5.3 of Golub and Van Loan~\cite{Golub89}).
When $A$ is Toeplitz (but not symmetric positive definite) we can
justify use of the semi-normal equations.
This is because the usual stable algorithms
for solving~(\ref{eq:Axb}) directly require $O(n^3)$ arithmetic operations,
but we can use the algorithm of \S\ref{sec:Main}
to compute (a numerical approximation $\cR$ to) $R$ in $O(n^2)$ operations,
and then solve the seminormal equations~(\ref{eq:seminormal}) in an additional
$O(n^2)$ operations.

From Theorem~\ref{thm:main}, we can compute
an upper triangular matrix $\cR$ such that~(\ref{eq:mainbound}) holds.
We can also compute an accurate approximation $\cd$ to
$d = A^Tb$ in $O(n^2)$ operations (using the obvious algorithm)
or in $O(n \log n)$ operations (using the Fast Fourier Transform).
Now solve the two triangular systems $\cR^Tw = d$ and $\cR x = w$.
We expect to obtain a result $\cx$ for which
\beq
\|\cx - x\|/\|x\| = O_n(\kappa^2\veps),			\label{eq:k2err}
\eeq
where $\kappa = \kappa(A)$, provided $\kappa^2\veps \ll 1$.
The residual $r = A\cx - b$ should satisfy
\beq
{\|r\| \over \|A\| \|x\|} = O_n(\kappa\veps),		\label{eq:k1err}
\eeq
because $\|A^Tr\| = \|A^TA\cx - A^Tb\| = O_n(\veps\|A^TA\|\|x\|)$.
From~(\ref{eq:k2err}), the method is {\em weakly stable} (according to
Definition~\ref{def:weak}), although we can not expect
the stronger bound~(\ref{eq:relerr}) to be satisfied.

The bounds~(\ref{eq:k2err}--\ref{eq:k1err}) are similar to those
usually given for the method of normal equations~\cite{Golub89},
not those usually given for the method of
semi-normal equations~\cite{Bjorck87,Paige73,Saunders72}. This is
because, in applications of the semi-normal equations, it is
usually assumed that $\cR$ is computed via an orthogonal
factorization of $A$, so there is a matrix $\hA$ such that
$\cR^T\cR = \hA^T\hA$ and
\beq
\|\hA - A\|/\|A\| = O_n(\veps).				\label{eq:k0}
\eeq
However, in
our case we only have $\|\cR^T\cR - R^TR\|/\|R^TR\| = O_n(\veps)$,
which implies the weaker bound
\beq
\|\hA - A\|/\|A\| = O_n(\kappa\veps)			\label{eq:k1}
\eeq
by Stewart's perturbation analysis~\cite{Stewart77,Stewart79}.

An alternative to the use of~(\ref{eq:seminormal}) was suggested
by Paige~\cite{Paige73}: compute $R$ such that\linebreak
$R^TR = AA^T$,
solve $R^TRw = b$ by solving two triangular systems,
then set $x = A^Tw$.	%
We prefer to use~(\ref{eq:seminormal}) because it is also
applicable in the rectangular (least squares) case~-- see
\S\ref{subsec:lsq}.

\subsection{Storage Requirements}
\label{subsec:Storage}

The algorithm described in \S\ref{sec:sneq}
for the solution of the semi-normal equations~(\ref{eq:seminormal})
requires working storage $O(n^2)$ words, because the upper triangular
matrix $R$ is not Toeplitz. However, it is possible to reduce
the storage requirement to $O(n)$ words. Recall that $\cR$ is
generated row by row. Thus, we can solve $\cR^Tw = d$ as $\cR$ is generated,
accumulating the necessary inner products with $O(n)$ storage.
We now have to solve $\cR x = w$ without having saved $\cR$. Provided
the $O(n)$ rotations defining the updates and downdates have been
saved, we can regenerate the rows of $\cR$ in reverse order
($n$, $n-1$, \ldots, $1$) and solve $\cR x = w$ as this is done.
A similar idea was used in~\cite{rpb079} to save storage in a systolic
implementation of the Bareiss algorithm.

The regenerated matrix $\ccR$ differs slightly from $\cR$ because
of rounding errors. Numerical experiments suggest (though we have
not proved) that $\|\ccR - \cR\| = O_n(\kappa\veps)$. 	%
If true, this would
imply that the method is weakly stable,
and~(\ref{eq:k2err}) would hold,
but the right hand side
of~(\ref{eq:k1err}) would need to be multiplied by $\kappa$.

An alternative is to use an idea which was suggested by
Griewank~\cite{Griewank92} in a different context.
Suppose we have a procedure $\calF(k)$ which generates $k$ consecutive
rows of $\cR$ in a forward direction (say rows $d+1,\ldots,d+k$),
and $\calB(k)$ which generates $k$
consecutive rows in a backward direction (say rows $d+k,\ldots,d+1$).
In each case $O(n)$
words of storage suffice for the initial conditions.
Then $\calB(2k)$ can be defined recursively~--
\begin{enumerate}
\item Generate rows $d+1,\ldots,d+k$ using $\calF(k)$,
	with row $d$ for initial conditions if $d > 0$,
	saving row $d+k$.
\item Generate rows $d+2k,\ldots,d+k+1$ using $\calB(k)$
	with row $d+k$ for initial conditions.
\item Generate rows $d+k,\ldots,d+1$ using $\calB(k)$
	with row $d$ for initial conditions if $d > 0$.
\end{enumerate}
From the discussion in \S\ref{sec:Main},
$\calF(k)$ requires $O(kn)$ operations and $O(n)$ storage.
Thus, we can prove by induction that
$\calB(k)$ requires $O(kn\log k)$ operations and $O(n\log k)$ storage.
Overall, to generate the $n$ rows of $\cR$ in reverse order
takes $O(n^2\log n)$ operations and $O(n\log n)$ storage.
Thus, at the cost of a factor $O(\log n)$ in the operation count,
we can reduce the storage requirements from $O(n^2)$ to $O(n \log n)$.
The numerical properties of this method are exactly the same as those
of the method which uses $O(n^2)$ storage, since exactly the same
rows of $\cR$ are computed.
The scheme described here is not optimal in its use of storage,
but is within a factor of two
of the optimal scheme described in~\cite{Griewank92}.

\subsection{Toeplitz Least Squares Problems}
\label{subsec:lsq}

If $A \in \Re^{m \times n}$ is Toeplitz with full rank $n$, then
the semi-normal equations~(\ref{eq:seminormal}) may be used to
solve the least squares problem
\beq
\min \|Ax - b\|_2.				\label{eq:Axb2}
\eeq
The use of semi-normal equations for the general full-rank
linear least squares problem
is discussed in detail by Bj\"orck~\cite{Bjorck87},
and the only significant difference in our case is that
the bound~(\ref{eq:k1}) holds instead of~(\ref{eq:k0}),
so an additional factor $\kappa$ appears in some of the terms in
the error bounds.

\subsection{Operation Counts}
\label{subsec:opcount}

We briefly estimate the number of arithmetic operations required by
the BBH algorithm and some of its competitors. We assume that the
Toeplitz matrix $A$ is real.
Some of the operation counts can be reduced by using fast Givens
transformations~\cite{Gentleman73,Golub89}, but for simplicity we
ignore this possibility.
We only count multiplications; the number of additions is comparable.

For simplicity, first consider the case $m = n$.
In the BBH algorithm, the computation of $R$ takes
$7n^2 + O(n)$ multiplications. The computation of $A^Tb$ takes
$n^2$ multiplications if done in the obvious way,
or $O(n\log n)$ multiplications if the FFT is used.
Solving two triangular systems takes $n^2 + O(n)$ multiplications.
Thus, to solve a Toeplitz linear system by the method of \S\ref{sec:sneq}
takes $9n^2 + O(n)$ multiplications, or $8n^2 + O(n\log n)$ if the FFT
is used. This is much cheaper than the
$19n^2 + O(n)$ %
multiplications of the
BBH $QR$ factorisation algorithm given in {\S}3 of~\cite{rpb092},
which computes $Q$ explicitly. Thus, considering both speed
and stability, it is best to avoid the computation of $Q$.

For the method of Nagy~\cite{Nagy93}, the multiplication count is
$16n^2 + O(n)$,
and for Cybenko's method~\cite{Cybenko87} it is $23n^2 + O(n)$.
The method TpH of~\cite{Gohberg94} requires $21n^2/2 + O(n\log n)$
real multiplications, and the method GKO of~\cite{Gohberg94} requires
$13n^2/2 + O(n\log n)$ complex multiplications.
Thus, the method of \S\ref{sec:sneq} should be faster
than any of these methods, although
GKO may be competitive if the Toeplitz matrix $A$ complex.

For the rectangular case $(m \ge n)$, the corresponding multiplication
counts (omitting low-order terms and assuming that the FFT is not used) are:

$2mn + 7n^2$ for the method of \S\ref{sec:sneq};

$2mn + 9n^2$ for the method of Nagy~\cite{Nagy93} using the
	semi-normal equations;	%

$2mn + 14n^2$ for the method of Nagy~\cite{Nagy93} using
	inverse factorisation (algorithm IF);

$9mn + 14n^2$ for the method of Cybenko~\cite{Cybenko87};

$13mn + 7n^2$ for the method of~\cite{rpb092},
	computing $Q$ explicitly.	%

\noindent For details of the components making up these operation counts,
see Table~4.1 of~\cite{Nagy93}.

\subsection{Iterative Refinement}
\label{subsec:refinement}

Iterative refinement (sometimes called iterative
improvement) can be used to improve an approximate solution $\cx$
to the linear system~(\ref{eq:Axb}) or the linear least
squares problem~(\ref{eq:Axb2}).
In practice this gives an accurate solution in a small number
of iterations so long as the residual is computed accurately
and the working precision is sufficient to ensure convergence.
For details
see~\cite{Bjorck87,Golub89,Golub66,Jankowski77,Skeel80,Wilkinson65}.

A related idea is to
improve the accuracy of $\cRb$ and $\cRt$ by using
Bj\"orck's ``Corrected Semi-Normal Equations''~\cite{Bjorck87,Park93}
or Foster's scheme of iterative improvement~\cite{Foster91}.
However, if the aim is simply to solve a linear system, then it is
more economical to apply iterative refinement directly to the system.

\subsection{Ill-conditioned Problems}
\label{subsec:illcond}

For very ill-conditioned Toeplitz linear systems and least squares problems,
it may be desirable to use regularisation, as discussed
in {\S}5 of~\cite{Nagy93}. We do not consider this here, except to note
that our algorithm can easily be modified to compute the Cholesky
factorisation of $A^TA + \alpha I$, where $\alpha$ is a
positive regularisation parameter.
Only small changes in equations~(\ref{eq:r11}--\ref{eq:Rt}) are required.

\section{Numerical Results}
\label{sec:Numer}

The algorithm described in \S\S\ref{sec:Main}--\ref{sec:sneq} has
been implemented in Pascal on an IBM~PC and DEC VAX.
In Table~\ref{tab:num1} we give some results for randomly chosen
$n \times n$ Toeplitz systems on a DEC VAX with
$\veps = 2^{-56} \simeq 1.4\times 10^{-17}$.
The elements $a_k$ defining the Toeplitz matrix
\[
A = \left(\begin{array}{ccc}
a_0	& \cdots	& a_{n-1}	\\
\vdots	& \ddots	& \vdots	\\
a_{1-n}	& \cdots	& a_0
\end{array}\right)
\]
were chosen from a normal distribution with specified
mean $\mu$ and standard deviation $\sigma$.
(The condition number $\kappa(A)$ tends to
increase with $|\mu/\sigma|$.)
The solution vector $x$ was chosen with normally distributed
components (mean~0) and the vector $b$ computed
from $b \leftarrow Ax$.
A consequence is that $\|x\|/\|b\|$ is unlikely
to be large, even if $A$ is poorly conditioned,
but this is typical of most applications.
Table~\ref{tab:num1} gives the condition number
\[ \kappa_1(R) = \|R\|_1\cdot\|R^{-1}\|_1,	\]
which is a rough approximation to $\kappa_2(R) = \kappa_2(A)$
(the 1-norm was used for computational convenience).

Table~\ref{tab:num1} gives
\begin{eqnarray*}
e_1 &=& {\|\cR^T\cR - A^TA\|_1 \over \veps\|A^TA\|_1},\\
e_2 &=& {\|\cx - x\|_2 \over \veps\kappa_1(R)^2\|x\|_2},\\
e_3 &=& {\|r\|_2 \over \veps\kappa_1(R)\|A\|_1\|x\|_2},
\end{eqnarray*}
where $r = A\cx - b$.
From Theorem~\ref{thm:main} and the bounds~(\ref{eq:k2err}--\ref{eq:k1err}),
we expect these quantities to be bounded by
(low degree) polynomials in $n$.
The results confirm this.

For comparison, the last column of Table~\ref{tab:num1}
gives $e_3^c$, the value of the normalised
residual $e_3$ obtained
via Cholesky factorization of $A^TA$.
It can be seen that $e_3$ is not much larger than $e_3^c$.

\begin{table}[htb]
\centerline{
\begin{tabular}{|c|l|l|l|l|l||l|} \hline
$n$	&$\mu/\sigma$
			&$\kappa_1(R)$
	& \hspace{0.5em}$e_1$
			& \hspace{0.5em} $e_2$
					& \hspace{0.5em} $e_3$
					     & \hspace{0.5em} $e_3^c$ \\ \hline
50	& 0.0		& 3.3{\E}3
	& 6.7		& 9.2{\E}-4	& 5.2{\E}-3	& 1.4{\E}-3\\
50	& 1.0		& 3.6{\E}3
	& 1.7{\E}1	& 6.2{\E}-3	& 2.7{\E}-2	& 2.1{\E}-3\\
50	& 1.0{\E}1	& 9.8{\E}2
	& 4.0{\E}1	& 3.1{\E}-1	& 3.7{\E}-1	& 9.4{\E}-2\\
50	& 1.0{\E}2	& 1.0{\E}4
	& 1.0{\E}2	& 2.0{\E}-1	& 5.3{\E}-1	& 7.0{\E}-2\\
50	& 1.0{\E}3	& 1.5{\E}5
	& 3.4{\E}1	& 4.6{\E}-1	& 3.0{\E}-1	& 4.6{\E}-2\\
50	& 1.0{\E}4	& 9.1{\E}5
	& 1.0{\E}2	& 1.0		& 1.2		& 9.2{\E}-2\\
50	& 1.0{\E}5	& 6.3{\E}6
	& 9.2		& 1.4{\E}-1	& 1.6{\E}-1	& 2.9{\E}-1\\ \hline
100	& 0.0		& 1.2{\E}3
	& 1.4{\E}1	& 6.2{\E}-4	& 4.7{\E}-3	& 1.6{\E}-3\\
100	& 1.0		& 1.8{\E}3
	& 1.5{\E}1	& 5.2{\E}-4	& 4.4{\E}-3	& 3.1{\E}-3\\
100	& 1.0{\E}1	& 5.2{\E}3
	& 8.7{\E}1	& 8.5{\E}-2	& 1.6{\E}-1	& 3.3{\E}-2\\
100	& 1.0{\E}2	& 3.3{\E}4
	& 1.2{\E}2	& 4.2{\E}-1	& 3.7{\E}-1	& 3.5{\E}-2\\
100	& 1.0{\E}3	& 4.0{\E}5
	& 1.4{\E}2	& 2.2{\E}-1	& 6.6{\E}-1	& 4.8{\E}-2\\
100	& 1.0{\E}4	& 8.8{\E}6
	& 1.5{\E}2	& 1.0		& 8.9{\E}-1	& 2.6{\E}-2\\
100	& 1.0{\E}5	& 7.8{\E}6
	& 3.8{\E}1	& 5.5{\E}-1	& 8.4{\E}-1	& 1.4{\E}-1\\ \hline
200	& 0.0		& 5.7{\E}3
	& 1.9{\E}1	& 1.3{\E}-4	& 1.7{\E}-3	& 8.7{\E}-4\\
200	& 1.0		& 1.2{\E}6
	& 8.4{\E}1	& 1.0{\E}-5	& 1.9{\E}-4	& 1.6{\E}-3\\
200	& 1.0{\E}1	& 5.6{\E}5
	& 1.6{\E}1	& 3.6{\E}-3	& 7.0{\E}-3	& 4.7{\E}-3\\
200	& 1.0{\E}2	& 2.4{\E}4
	& 2.1{\E}2	& 3.0		& 2.7		& 1.4{\E}-1\\
200	& 1.0{\E}3	& 7.9{\E}5
	& 1.4{\E}1	& 8.2{\E}-2	& 6.6{\E}-2	& 5.2{\E}-2\\
200	& 1.0{\E}4	& 5.5{\E}6
	& 2.8{\E}2	& 6.1{\E}-1	& 5.0{\E}-1	& 4.3{\E}-2\\
200	& 1.0{\E}5 	& 1.3{\E}8
	& 3.6{\E}2 	& 1.8{\E}-1	& 3.8{\E}-1	& 4.3{\E}-2\\ \hline
\end{tabular}
}
\caption{Weakly stable solution of Toeplitz systems}
\label{tab:num1}
\end{table}

We also tried Toeplitz matrices $A$ with some singular principal
submatrices (e.g.~$A$ with $a_{-1} = a_0 = a_1$). The results were
similar to those given in Table~\ref{tab:num1}.

\section{Conclusion}
\label{sec:Conc}

The method described in \S\ref{sec:sneq}
for the solution of general nonsingular Toeplitz or Hankel linear
systems\footnote{Similar remarks apply for full-rank
Toeplitz or Hankel least squares problems.}
requires $O(n^2)$ operations,
is weakly stable, and makes no assumption
about the conditioning of submatrices of $A$.
We do not know any other methods which have been proved to be
stable or weakly stable and have worst-case time bound $O(n^2)$.
Algorithms which involve pivoting and/or
look-ahead~\cite{Chan92a,Freund93d,Gutknecht93a,Hansen93,Sweet93}
may work well in practice,
but seem to require worst-case overhead $O(n^3)$ to ensure stability.

Our method should be faster than $O(n^3)$ methods which ignore the
Toeplitz structure,
even if the working precision has to be increased (i.e.~$\veps$
reduced) in our method to ensure that $\kappa^2\veps \ll 1$.
Storage requirements are $O(n^2)$, but may be reduced to
$O(n \log n)$ or $O(n)$ as discussed in \S\ref{subsec:Storage}.

We have not looked in detail at all the fast Toeplitz $QR$ factorization
algorithms mentioned in \S\ref{sec:Intro}. It is quite likely that
some of these give weakly stable algorithms for
the computation of the Cholesky factor $R$ of $A^TA$,
but no proofs have been published.

It remains an open question whether there is a fast Toeplitz $QR$
algorithm which is backward stable for the computation of $Q$
and $R$ (or $R^{-1}$).
Such an algorithm would give a stable algorithm for the solution of
$Ax = b$ without recourse to the semi-normal equations.

\subsection*{Acknowledgements}

Thanks to
Greg Ammar,
George Cybenko,
Lars Eld\'en,
Andreas Griewank,
Georg Heinig,
Franklin Luk,
Haesun Park,
and Douglas Sweet,
for helpful comments on earlier versions of this paper.
James Bunch,
Roland Freund,
Martin Gutknecht,
Vadim Olshevsky
and James Varah
kindly sent us relevant reprints and preprints.

\bigskip\noindent The literature on Toeplitz and related linear systems
is voluminous, so the following bibliography is by no means complete.
We have attempted to include a representative sample of the
recent literature and some important historical references.

{\small			%

}%

\bigskip

\section*{Postscript}

At the time of writing this paper we were not aware of a paper by S. Cabay
and R. Meleshko [A weakly stable algorithm for Pad\'e approximants and the
inversion of Hankel matrices, {\em SIAM J. Matrix Anal. Appl.} 15 (1993),
735--765], which gives a different weakly stable algorithm for the inversion
of Toeplitz/Hankel matrices, usually (but not always) in $O(n^2)$ operations.

\end{document}